\documentclass[12pt]{article}%
\usepackage{amsfonts}
\usepackage{amssymb}
\usepackage{amsmath}%
\usepackage{graphicx}
\providecommand{\U}[1]{\protect\rule{.1in}{.1in}}
\newtheorem{theorem}{Theorem}

\newcommand{\M}{\mathcal{M}}

\begin{document}

\title{The exponents of $p$-groups of maximal class and their Schur multipliers}
\author{Michael Vaughan-Lee}
\date{July 2026 }
\maketitle

\section{Introduction}

In this note we prove the following two theorems.

\begin{theorem}
Let $P$ be a $p$-group of maximal class of order $p^n$ where $p,n>2$. If $n>p+1$ then the exponent of $P$ is
$p^k$, where $k$ is the least integer greater or equal to $(n-1)/(p-1)$. If $n=p+1$ then $P$ has exponent $p^2$.
And if $n\leq p$ then $P$ has exponent $p$ or $p^2$, with both possibilities arising.
\end{theorem}

\begin{theorem}
Let $P$ be a $p$-group of maximal class of order $p^n$ where $p,n>2$, and let $\M(P)$ be the Schur multiplier
of $P$. Then $\M(P)$ has exponent at most $p^k$, where $k$ is the least integer greater than or equal to
$(n-2)/(p-1)$.
\end{theorem}

Note that if we let
\[
P=\langle a,b\,|\,a^{p}=1,\,(ab)^{p}=1,\,[b^{a^{i}},b]=1\;(1\leq i<p),\,[b,_{n-1}a]=1\rangle
\]
then $P$ is a group of maximal class of order $p^{n}$. Its Schur multiplier has exponent $p^k$
where $k$ is the least integer greater than or equal to $(n-2)/(p-1)$. See \cite[Theorem 2]{hadi}.

These results are undoubtedly well known, and in fact they follow fairly easily from the detailed results
on $p$-groups of maximal class to be found in Chapter 3 of Charles Leedham-Green
and Sue McKay's book \cite{leedhamgm}. However these results do not appear in \cite{leedhamgm},
and nor do they appear in Blackburn's seminal paper \cite{blackburn} on groups of maximal class.

\section{Preliminaries}

Let $P$ be a $p$-group of maximal class and order $p^{n}$. The nilpotency
class of $P$ is $n-1$. Let%
\[
P>P_{2}>P_{3}>...>P_{n-1}>P_{n}=\{1\}
\]
be the lower central series of $P$. For $i=2,3,\ldots,n-2$ we let the $2$-step
centralizer $K_{i}$ be the centralizer in $P$ of $P_{i}/P_{i+2}$. The
subgroups $K_{2}$, $K_{3}$, \ldots, $K_{n-2}$ all contain the derived subgroup
$P_{2}$ and are all maximal subgroups of $P$. It turns out that $K_{2}%
=K_{3}=...=K_{n-3}$. There is a proof of this fact in Charles Leedham-Green
and Sue McKay's book \cite{leedhamgm}. The result is not explicitly stated in
their book, but it follows easily from Corollary 3.2.7, Theorem 3.2.11, and
Theorem 3.3.5. So the union of all the $2$-step centralizers is just
$K_{2}\cup K_{n-2}$, and we can pick an element $s\in P$ which lies outside
this union. We also pick an element $s_{1}$ in $K_{2}\backslash P_{2}$. So $P$
is generated by $s$ and $s_{1}$.

We let $P_{1}=\langle s_{1}\rangle P_{2}$, so that%
\[
P>P_{1}>P_{2}>P_{3}>...>P_{n-1}>P_{n}=\{1\}.
\]
For $i=1,2,\ldots,n-2$ we let $s_{i+1}=[s_{i},s]$, so that $s_{i}$ generates
$P_{i}$ modulo $P_{i+1}$ for $i=1,2,\ldots,n-1$.

\section{Proof of Theorem 1}

First consider the case when $n>p+1$. Then, by \cite[Theorem 3.3.5]{leedhamgm},
$P$ has positive degree of commutativity, which implies that all the 2-step
centralizers are equal. So, by \cite[Lemma 3.3.1]{leedhamgm}, if $t \notin P_1$
then $t^{p^2}=1$. On the other hand, \cite[Corollary 3.3.6]{leedhamgm} states
that if $1\leq i \leq n-p+1$ then $P_i^p=P_{i+p-1}$. It follows that $P_1^{p^k}=1$
if and only if $1+k(p-1)\ge n$, and this proves Theorem 1 in the case when $n>p+1$.

The case when $n=p+1$ is slightly trickier. By \cite[Proposition 3.3.2]{leedhamgm} we see that if
$n=p+1$ then $P/P_{n-1}$ has exponent $p$, and so $P$ has exponent at most $p^2$. We show that
$P$ has exponent exactly $p^2$ by showing that $[s_1,_{p-1}s]$ is a product of $p$-th powers.

Note that if $2\leq i \leq n-3$ then $s_1 \in K_i$, so that $[s_i,s_1]\in P_{i+2}$. 
Putting this in another way, we see that if $a_1,a_2,\ldots ,a_{i+1}$ are elements of the set 
$\{s,s_1\}$, and if two or more of $a_1,a_2,\ldots ,a_{i+1}$ are equal to $s_1$ then 
\[
[a_1,a_2,\ldots ,a_{i+1}]\in P_{i+2}.
\]
Since $P$ is nilpotent of class $p$ this implies that any commutator in $P$ of
weight $p$ is trivial if three or more of the entries in the commutator are
equal to $s_1$.

Now let $F$ be the free group of rank two generated by $s,s_1$, and apply the
Hall collection process to $(ss_1)^p$. Working modulo $\gamma _{p+1}(F)$
we obtain
\[
(ss_1)^p=s^ps_1^p[s_1,s]^{p(p-1)/2}c_4^{\alpha _4}c_5^{\alpha _5}\ldots c_m^{\alpha _m}
\]
where $c_4,c_5,\ldots ,c_m$ are basic commutators in $s,s_1$ of weight at most $p$.
If $c_i$ has weight $w$ then $\alpha _i$ is an integral linear combination of the
binomial coefficients $\binom{p}{1},\binom{p}{2},\ldots ,\binom{p}{w}$.
Note that the basic commutator $[s_1,_{p-1}s]$ occurs in this product with exponent 1,
and that the basic commutator $[s_1,_{p-2}s,s_1]$ occurs with exponent $p^2-p-1$.

Let $N$ be the subgroup of $P$ generated by $p$-th powers in $P$.
Evaluating this equation in $P$ we use the fact
that if $c_i$ has weight less than $p$ then $\alpha _i$ is divisible by $p$.
We also use that fact that if $c_i$ has weight $p$ and three or more entries
$s_1$ then $c_i$ exaluates to 1 in $P$. So we obtain
\[
[s_1,_{p-1}s][s_1,_{p-2}s,s_1]^{-1}c \in N
\]
where $c$ is a product of powers of basic commutors
$[[s_1,_is],[s_1,_js]]$ where $i>j>0$ and $i+j=p-2$. But in $P$
\[
[[s_1,_is],[s_1,_js]]=[s_1,_{p-2}s,s_1]^{(-1)^j},
\]
So we see that applying Hall collection to $(ss_1)^p$ gives us an equation
\[
[s_1,_{p-1}s][s_1,_{p-2}s,s_1]^k \in N
\]
for some integer $k$. (Working modulo $p$ it seems that $k=(p+1)/2$, but I have no proof of this.)
Note that this equation does not depend on the particular choice for $s_1$, so substituting
$s_1^2$ for $s_1$ in the equation we obtain
\[
[s_1^2,_{p-1}s][s_1^2,_{p-2}s,s_1^2]^k \in N
\]
However these two equations imply that $[s_1,_{p-1}s] \in N$, so that $P$ has exponent $p^2$.

Finally consider the case when $n<p+1$. By \cite[Proposition 3.3.2]{leedhamgm}, 
$P/P_{n-1}$ has exponent $p$, and so $P$ has exponent at most $p^2$. The
group
\[
\langle a,b\,|\,a^{p}=1,\,(ab)^{p}=1,\,[b^{a^{i}},b]=1\;(1\leq i<p),\,[b,_{n-1}a]=1\rangle
\]
is  group of maximal class of order $p^n$ and it has exponent $p$ since it is generated
by elements of order $p$, and has class $n-1<p$.
And the group
\[
P=\langle a,b\,|\,a^{p}=[b,_{n-2}a],\,(ab)^{p}=1,\,[b^{a^{i}},b]=1\;(1\leq i<p),\,[b,_{n-1}a]=1\rangle
\]
is a group of maximal class of order $p^n$ with exponent $p^2$.

\section{Proof of Theorem 2}

Let $P$ be a group of maximal class of order $p^n$ where $p,n>2$, and let $P=F/R$ where
$F$ is the free group of rank 2 generated by $s,s_1$. The Schur cover of $P$ is $F/[R,F]$
and $\M(P) \leq \gamma_3(F/[R,F])$.  We prove Theorem 2 by showing that $\gamma_3(F/[R,F])$
has exponent at most $p^k$ where $k$ is the least integer greater than or equal to
$(n-2)/(p-1)$.

First consider the case when $n\leq p+1$. Then $P$ has class at most $p$, and $F/[R,F]$ has
class at most $p+1$. We will show that $\M(P)$ has exponent $p$. Since $\M(P)\leq \gamma_3(F/[R,F])$,
it is enough to show that $\gamma_3(F)^p \leq [R,F]$. And since $\gamma_3(F/[R,F)$ has
class less than $p$ it is enough to show that if $m\ge 3$ and $a_1,a_2,\ldots ,a_m\in F$
then $[a_1,a_2,\ldots ,a_m]^p \in [R,F]$.

So let $m\ge 3$ and let $a_1,a_2,\ldots ,a_m\in F$. By \cite[Proposition 3.3.2]{leedhamgm},
$a_m^p \in \gamma_{n-1}(F)R$, and so (since $\gamma_{n+1}(F) \leq [R,F]$)
\[
[a_1,a_2,\ldots ,a_{m-1},a_m^p] \in [R,F].
\]
Now
\[
[a_1,a_2,\ldots ,a_{m-1},a_m^p]=[a_1,a_2,\ldots ,a_m]^p\text{ modulo }\gamma_{m+1}(F)^p\gamma_{m-1+p}(F),
\]
and $\gamma_{m-1+p}(F) \leq [R,F]$. So $\gamma_m(F)^p \leq \gamma_{m+1}(F)^p$ modulo $[R,F]$,
and since this true for all $m\ge 3$ we concude that $\gamma_3(F)^p \leq [R,F]$.

Next, suppose that $n>p+1$. By \cite[Corollary 3.3.6]{leedhamgm}, $P_i^p=P_{i+p-1}$ for 
$1\leq i\leq n-p+1$. This implies that $\gamma_i(F)^p\leq \gamma_{i+p-1}(F)R$ for
$2\leq i\leq n-p+1$. So
\[
[\gamma_i(F)^p,F]\leq [\gamma_{i+p-1}(F),F][R,F]=\gamma_{i+p}(F)[R,F].
\]
Let $x\in \gamma_i(F)$ and let $y\in F$. Then, working modulo $\gamma_{i+p}(F)$,
$[x^p,y]=[x,y]^pz$ where $z$ is a product of $p$-th powers of elements in
$\gamma_{i+2}(F)$. So
\[
\gamma_{i+1}(F)^p\leq \gamma_{i+2}(F)^p\gamma_{i+p}(F)[R,F].
\]
Since this equation holds true for all $2\leq i \leq n-p+1$, if follows that 
\[
\gamma_{i+1}(F)^p\leq \gamma_{i+p}(F)[R,F]
\]
for all $i\ge 2$. So, working modulo $[R,F]$, 
\[
\gamma_3(F)^p \leq \gamma_{p+2}(F),\;\gamma_3(F)^{p^2} \leq \gamma_{2p+1}(F),\;\gamma_3(F)^{p^3} \leq \gamma_{3p}(F),
\]
and so on. In general if $k\ge 1$ then $\gamma_3(F)^{p^k} \leq \gamma_{3+k(p-1)}(F)$
modulo $[R,F]$. Since $\gamma_{n+1}(F)\leq [R,F]$ this implies that if $k\ge (n-2)/(p-1)$
then $\gamma_3(F)^{p^k} \leq [R,F]$.
This completes the proof of Theorem 2, since $\M(P)\leq \gamma_3(F/[R,F])$.

\end{document}